\documentclass[12pt]{article}
\usepackage{amssymb,latexsym,theorem}
\newcommand{\Z}{\mathbb Z}
\newcommand{\R}{\mathbb R}
\newcommand{\isodim}{\mathit{isodim}}
\newcommand{\GL}{\mathit{GL}}
\newcommand{\rank}{\mathit{rank}}
\newcommand{\colone}{\mathit{col}1}
\newcommand{\gram}{\mathit{gram}}
\newcommand{\kmax}{\mathit{kmax}}
\newcommand{\mix}{\mathrm{mix}}
\newcommand{\priso}{\mathrm{pr}_\mathrm{iso}}
\newcommand{\transpose}{\mathrm{transpose}}
\newcommand{\trickle}{\texttt{\textup {trickledown}}}
\newcommand{\diso}{\mathit{diso}}
\newcommand{\qed}{\unskip\nobreak\hfill\hbox{ $\Box$}}

\newtheorem{Theorem}[subsection]{Theorem}
\newtheorem{Lemma}[subsection]{Lemma}

\theorembodyfont{\normalfont}
\newtheorem{Remark}[subsection]{Remark}

\begin{document}
\sloppy
\title{Complexity of the Havas, Majewski, Matthews LLL Hermite Normal Form
algorithm}
\author{Wilberd van der Kallen}
\maketitle
\section{Summary}
We consider the complexity of the LLL HNF algorithm (\cite{HMM}).
This algorithm takes as input an $m$ by $n$ matrix $G$ of integers
and produces as output a matrix $b\in \GL_m(\Z)$ so that $A=bG$
is in
Hermite
normal form (upside down).
The analysis is similar to that of an extended LLL algorithm as given in
\cite{vdK}.
\section{Sketch of the argument}
Let $B$ be the maximum of the entries of $GG^\transpose$.
The main issue is whether we can estimate the entries of $b$
in terms of $B$, $m$, $n$ during the algorithm.
The entries of $A$ can then be estimated through $A=bG$.
As they do not affect $b$, we may remove from $G$ all columns that do not
contribute a pivot to  the
Hermite normal form.
Once that is done, $A$ has as many columns as its rank, and at the end
of the algorithm the product
of its pivots is the covolume of the lattice spanned by its rows.
This covolume can be estimated in terms of a minor of $G$,
which by Hadamard is at most $ B^{\rank/2}$.
This replaces the estimate $d_i\leq B^i$ of \cite{LLL}.
As in \cite{vdK} we use the ordinary Euclidean inner product
$(\ ,\ )$ for the rows of $b$, but also an inner product
$\langle v,w \rangle=(Gv,Gw)$.
(For the new $G$ which has a rank equal to its number of columns.)
The vectors $v$ with $\langle v,v\rangle=0$ are called isotropic.
We let $\priso$ be the orthogonal projection
according to $(\ ,\ )$ of $\R^m$
onto
the isotropic subspace and put
$$(v,w)_\mix=(\priso v,\priso w)+\langle v ,w \rangle .$$
One estimates the ratio between  $(v,v)$ and $(v,v)_\mix$
and one estimates $(b_i,b_i)_\mix$ for any row $b_i$ of $b$.
The required estimate of the entries of $b$
follows from this, at least at the
end of the algorithm.
The algorithm computes a Hermite normal form first
for the top $\kmax$ rows of $G$, starting with $\kmax=1$,
 and increases $\kmax$ in steps of one.
Each time just before
one wants to increase $\kmax$ the analysis sketched above applies.
Right after one  wants to increase $\kmax$
 we enter a stage which we will emulate with a procedure
called $\trickle$, which we analyze as in \cite{vdK}. It is followed
by an ordinary LLL stage and then we get back to increasing $\kmax$.
Therefore no accidents happen.

\section{The analogy with an extended LLL algorithm}
Let $e_1$,\ldots,$e_m$ be the standard basis of $\R^m$.
The  Gram matrix
$\gram=(\langle e_i,e_j\rangle )_{i,j=1}^m$
belongs to a positive semidefinite inner product $\langle\ ,\ \rangle$ on
$\R^m$.
Note that $\gram$ has integer entries.
Let $\rank $ be the rank of $G$ and assume that we have removed from $G$
the columns that do not contribute a pivot.
(In this paper a pivot is an entry of $A$ that is
the first nonzero entry in its row and also in its column.)
Put   $\isodim=m-\rank$ and let $b_i^*$ denote the $i$-th Gram-Schmidt
vector in the following sense. We have $b_i^*\in (b_i+\sum_{j=1}^{i-1}\R b_i)$
and if  $1\leq j<i$, $j\leq \isodim$ then
$(b_i^*, b_j)=0$, but if $1\leq j<i$, $j> \isodim$ then
$\langle b_i^*, b_j\rangle =0$.
With those notations the output satisfies:
\begin{enumerate}
\item The first $\isodim$ rows $b_i$ of $b$ are isotropic.
\item With respect to $(\ ,\ )$ the first $\isodim$ rows of $b$ form an LLL
reduced basis of $\sum_{j=1}^{\isodim}\Z b_i$.
\item  The last $\rank$ rows of $b$ form a basis of the lattice they span,
and this lattice contains no
nonzero isotropic vector. Furthermore,
$|\langle b_i^*,b_j\rangle |\leq \langle b_i^*,b_i\rangle $
for $\isodim+1\leq i<j\leq m$.
\item For $1\leq i\leq \isodim$, $i<j\leq m$ we have
$|(b_i^*,b_j)|\leq 1/2(b_i^*,b_i)$.
\end{enumerate}

\begin{Remark}
These properties are very much in the spirit of \cite{LLL} and will thus allow
us to estimate the $(b_i,b_i)_\mix$ in a traditional manner.
Of course the output satisfies more properties, as
actually $bG$ has Hermite normal form, but we suppress that now.
\end{Remark}

\section{Stages}
Now that we have a way to look at the final result, let us discuss how
we view things along the way.
The moment that $k$ wants to go beyond $\kmax$ is special.
(As in \cite{C} we use $\kmax$ to denote the maximum value that $k$ has
attained.)
At this moment one can estimate everything in the same manner as at the end.
It is followed by a stage which we emulate by the procedure \trickle,
and which ends
when
 a new pivot or a new isotropic
vector appears in the actual Havas, Majewski, Matthews LLL Hermite Normal Form
algorithm.
We then have an estimate of $b$ as in \cite{vdK}.
After this one basically runs an ordinary LLL algorithm for the
inner product $(\ ,\ )_\mix$ until $k$ wants to go beyond $\kmax$ again.

What makes it all rather technical is that $(\ ,\ )_\mix$ always depends
on which pivots and which independent isotropic vectors have
 been found. During \trickle\ one needs to take into account that one
is dealing with an MLLL in the sense of \cite{P}, with the added
complication that it is an \emph{extended} MLLL
algorithm in that one also requires
the transformation matrix $b$. It is the latter which makes that one can
not refer to \cite{P} for the analysis.

\section{Running LLL}
In the definition of  the inner product $\langle\ ,\ \rangle$ we work with a
$G$ from which all columns have been removed  where no pivot
has been found yet.
In the ordinary LLL stages we will have,
\begin{itemize}
\item An
integer matrix $b$ of
determinant one,
\item Integers
$k$, $\kmax$, $1\leq k\leq \kmax\leq m$,
\item An integer $\isodim\geq0$,  so that
the first $\isodim$ rows $b_i$ of $b$ span the isotropic subspace of
$\sum_{j=1}^{\kmax}\R b_j$.
\end{itemize}
(Initialize with $k=\kmax=1$ and $\isodim=0$.)

Let $\priso$ be the orthogonal projection according to $(\ ,\ )$ of $\R^m$
onto
$\sum_{j=1}^{\isodim}\R b_j$ and put
$$(v,w)_\mix=(\priso v,\priso w)+\langle v ,w \rangle .$$
Let $\mu_{i,j}$ be defined for $i>j$ so that
$$b_i=b_i^*+\sum_{j=1}^{i-1}\mu_{i,j} b_j^*.$$

The first standard assumption is then that, with respect to
$(\ ,\ )_\mix$, the first
$k-1$ rows of $b$ form an LLL reduced basis of $\sum_{j=1}^{k-1}\Z b_j$,
except that one does \emph{not} require
$$|b_{i}^*+\mu_{i,i-1}b_{i-1}^*|_\mix^2
\geq 3/4 |b_{i-1}^*|_\mix^2$$ when $i>\isodim$,
and that the usual condition  $|\mu_{i,j}|\leq 1/2$ is weakened to
 $|\mu_{i,j}|\leq 1$ for $j>\isodim$.
And the second standard assumption is that, as in \cite{C}, the first
$\kmax$ rows of $b$ form a basis of  $\sum_{j=1}^{\kmax}\Z e_i$.

We run the LLL algorithm with respect to $(\ ,\ )_\mix$,
except that one leaves out many swaps.
(From now on we suppress mentioning the annoying weakening of the condition on
the $\mu_{i,j}$.)
Leaving out swaps is harmless, as the size estimates in \cite{LLL}
for the $\mu_{ij}$ \textit{etcetera}
  do not require that one executes a swap whenever such
is recommended by the LLL test.
Running LLL with   respect to $(\ ,\ )_\mix$
roughly amounts to running two LLL algorithms, one for $(\ ,\ )$
and one for $\langle\ ,\ \rangle$.
That is why the pseudo-code in \cite{HMM} makes the distinction between
$\colone=n+1$
and $\colone\leq n$.
One runs LLL until $k$ tries to go to $\kmax+1$.
If $\kmax=m$ we are through.
If $\kmax< m$
then one should realize that because of the removal of columns from $G$
the row $e_{\kmax+1}G$ will be dependent on the earlier ones.
So we enter an extended MLLL, which we emulate with $\trickle$.

\section{Estimates}
We want to give estimates by changing \cite{vdK} minimally.
Thus let $B\geq2$ so that the entries of $\gram$ are at most $B$.
Our main result is
\begin{Theorem}
All through the algorithm all entries have bit length
${\cal O}(m \log (m B))$.
\end{Theorem}
We do not care about the constants in this estimate.
We leave to the reader the easy task of estimating the number of
operations on the entries in the manner of \cite{LLL}.
One finds that ${\cal O}((m+n)^4 \log (m B))$ such operations will do.

\subsection{Determinants}
Let $\gram_\mix$ be
the Gram matrix $((e_i,e_j)_\mix)$ with respect to
$e_1$, \ldots, $e_\kmax$.
Its entries are at most ${B+1}$.
 With Hadamard this gives
$$|\det(\gram_\mix)|\leq (\sqrt m (B+1))^{m}$$
and the same estimate holds for its subdeterminants.
We claim that the determinant of $\gram_\mix$ is an integer, so that we
also get this upper bound for the entries of $\gram_\mix^{-1}$.
To see the claim, consider as in \cite{P}
the inner product $(\ ,\ )_\epsilon$
given by $(v,w)_\epsilon=\epsilon(v,w)+\langle v,w\rangle$.
Its Gram matrix has a determinant which is a polynomial
$\det_\epsilon$ of $\epsilon$  with integer coefficients.
One may also compute $\det_\epsilon$
with respect to a basis which is obtained from
$e_1$, \ldots, $e_\kmax$ through an orthogonal transformation
matrix. By diagonalizing the Gram matrix of
$\langle\ ,\ \rangle$ we see that $\det(\gram_\mix)$ is the coefficient
of $\epsilon^\isodim$ in $\det_\epsilon$.\qed

\begin{Lemma}\label{mix2euc}
For $v\in\R^m$ one has $$(v,v)_\mix\leq m(B+1)(v,v)$$ and for
$v\in \sum_{j=1}^{\kmax}\R e_j$ one has
$$(v,v)\leq m(\sqrt m (B+1))^{m}(v,v)_\mix.$$
\end{Lemma}
\subsubsection*{Proof}
The supremum of $\{\,(v,v)_\mix\mid(v,v)=1\,\}$ is the
largest eigenvalue of the gram matrix of $(\ ,\ )_\mix$ with respect to
$e_1,\ldots,e_m$.  The
largest eigenvalue is no larger than the trace of this matrix. So it
is at most $m(B+1)$. Similarly the largest eigenvalue of $\gram_\mix^{-1}$
it is at most $m(\sqrt m (B+1))^{m}$.
\qed

\subsection{Vectors}
Now put $$\diso_i=\prod_{j=1}^i (b_j^*,b_j^*)$$ for $i\leq \isodim$
and $$d_i=\prod_{j=1}^i \langle b_{j+\isodim}^*,b_{j+\isodim}^*\rangle$$
for $i\leq \rank$.
As far as $d_i$ is concerned we may compute modulo isotropic vectors,
or also with $(\ ,\ )_\mix$.
Indeed $$\langle b_{j+\isodim}^*,b_{j+\isodim}^*\rangle=
(b_{j+\isodim}^*,b_{j+\isodim}^*)_\mix$$ for $1\leq j\leq \rank$.
Both $\diso_i$ and $d_j$ are integers and they descend when applying LLL.
(Throughout we assume familiarity with the arguments in \cite{LLL}.)
In fact the $\langle b_{j+\isodim}^*,b_{j+\isodim}^*\rangle$ are themselves
squares of integers. (Squares of the pivots of the moment.)

One may also
compute $\det(\gram_\mix)$ with the $b_i^*$ basis,
as the transition matrix has determinant one. From
that one sees that it is just $\diso_\isodim d_\rank$.
So we get
$\diso_\isodim\leq (\sqrt m (B+1))^{m}$. In fact, for $i\leq\isodim$
one has the same estimate
$$\diso_i\leq (\sqrt m (B+1))^{m}$$
because $i$ was equal to $\isodim$ earlier in the algorithm and LLL only
makes $\diso_i$ go down.
We have that $$d_i\leq B^\rank$$ because this is so when a pivot
is created and LLL only makes it descend. (Recall that the pivots are integers
whose product is at most $B^{\rank/2}$, while $d_i$ is a product of some
squared pivots
$\langle b_{j+\isodim}^*,b_{j+\isodim}^*\rangle$.
The
\trickle\ part also makes pivots descend.)

\begin{Lemma}\label{mixsmall}
Let $1\leq i\leq \kmax$. Then
 $$(\sqrt m (B+1))^{-m}\leq(b_{i}^*,b_{i}^*)_\mix\leq (\sqrt m (B+1))^{m}$$
and if $C\geq1$ is such that
$|\mu_{ij}|^2\leq C$ for $1\leq j<i$ then
$$(b_{i},b_{i})_\mix\leq mC(\sqrt m (B+1))^{m}$$
\end{Lemma}

\subsection*{Proof}
Use the estimates of $\diso_i$, $d_i$.
\qed

\subsection{Preserved estimates}
Put $C=(4mB)^{5m}$.
\begin{Lemma}\label{preserved}
The following estimates hold between applications of
\trickle\ (each time $k$ changes)
\begin{enumerate}
\item $\diso_i\leq (\sqrt m (B+1))^{m}$ for $i\leq\isodim$,
\item $d_i\leq B^\rank$ for $i\leq\rank$,
\item $(b_{i},b_{i})_\mix\leq mC(\sqrt m (B+1))^{m}$ for $i\neq k$,
\item $(b_{k},b_{k})_\mix\leq m^2 9^mC(\sqrt m (B+1))^{3m}$,
\item $|\mu_{i,j}|\leq 1$ for $1\leq j<i<k$,
\item $|\mu_{k,j}|\leq 3^{m-k}\sqrt {mC}(\sqrt m (B+1))^{m}$ for $1\leq j<k$,
\item $|\mu_{i,j}|\leq \sqrt {mC}(\sqrt m (B+1))^{m}$ for $1\leq j<i>k$.
\end{enumerate}
\end{Lemma}

\subsection*{Proof}
That these are preserved under LLL follows as in \cite{LLL}, so one
has to check that they hold right after \trickle. Given our earlier estimates,
this will be straightforward once we have shown that,
at that moment, $|\mu_{ij}|^2\leq C$.
Note that one could insert steps in the
algorithm to reduce to the case $C=1$ instead of the outrageously pessimistic
$C=(4mB)^{5m}$.
\qed

\section{Description of \trickle}

Before we can do estimates concerning \trickle\ we must describe it.
One starts with having
$k=\kmax+1\leq m$.
(So we look at the moment that $\kmax$ should be increased.)
 Consider the lattice
generated by $b_1,\ldots,b_{\kmax+1}$ where  $b_{\kmax+1}=e_{\kmax+1}$.
As $e_{\kmax+1}G$ is dependent on the earlier rows of $G$ now,
 this lattice contains a
nonzero vector $v$ with $(v,v)_\mix=0$.
Modulo $\R v$ the vector $b_{k}$
 is linearly dependent on the $b_i$ with $i<k$.
Changing the basis of $\Z b_{k-1}+\Z b_{k}$ we can achieve that
modulo $\R v$ the vector $b_{k-1}$
 is linearly dependent on the $b_i$ with $i<k-1$.
Then lower $k$ by one and repeat until $k=\isodim+1$,
where $\isodim$ is the one from before the present \trickle.
Or stop when the
Havas, Majewski, Matthews LLL Hermite Normal Form algorithm
produces a new pivot.
 If
$b_k=b_{\isodim+1}$ is itself isotropic we increase
$\isodim$ by one and pass to
a new $(\ ,\ )_\mix$.
If a new pivot has been created we add back the relevant column to $G$ and
again pass to
a new $(\ ,\ )_\mix$.
This describes \trickle.

One may worry about the fact that \trickle\ does not trace the
Havas, Majewski, Matthews LLL Hermite Normal Form algorithm
faithfully.
We are close enough though. (And our replacement is has worse estimates than
the original.)
We are just leaving out some size reductions
and we are taking together some swaps and reductions
 that make up the required change of
basis of $\Z b_{k-1}+\Z b_{k}$. The change of basis is the one coming from
an extended euclidean algorithm.
Thus we will further ignore that \trickle, which we took from \cite{vdK},
does not quite trace  this stage of the
Havas, Majewski, Matthews algorithm. We simply blame their algorithm.

\section{Estimates during \trickle}
We look in more detail.
Upon entering \trickle\
we freeze the old $\isodim$, $\kmax$ and the $b_i^*$, even though
the $b_i$ will change. We also do not change $(\ ,\ )_\mix$.
Let $\mu_{i,0}$ stand for $(e_{\kmax+1},b_i)$ and let
$\mu_{i,j}$ stand for
$(b_j^*,b_i)_\mix/(b_j^*,b_j^*)_\mix$ if $j>0$.
Note that initially $|\mu_{i,j}|\leq 1$ for $i\leq \kmax$,
$0\leq j\leq\kmax$. We will estimate  $|\mu_{i,j}|$ as $k$ descends.
The key point is that we can also estimate $\mu_{i,0}$. This compensates for
the fact that $(\ ,\ )_\mix$ is degenerate
on $\sum_{i=1}^{\kmax+1}\R e_i$. By combining $\mu_{i,0}$
with $(\ ,\ )_\mix$ we will be able to estimate $(b_i,b_i)$.
It is to explain the estimate of $\mu_{i,0}$ that we prefer to work 
with \trickle. 

Say $k>\isodim+1$ and modulo $\R v$ the vector $b_{k}$
 is linearly dependent on the $b_i$ with $i<k$.
Let us compute with $b_k$, $b_{k-1}$ modulo $V=\R v+\sum_{i=1}^{k-2}\R b_i$.
We have $b_k\equiv \mu_{k,k-1}b_{k-1}^*$ and $b_{k-1}\equiv b_{k-1}^*$
modulo $V$.
With the extended euclidean algorithm of \cite{C}
we find an integer matrix
$\pmatrix{\alpha&\beta\cr\gamma&\delta}$ of determinant one so that
$\pmatrix{\alpha&\beta\cr\gamma&\delta}\pmatrix{1\cr \mu_{k,k-1}}=
\pmatrix{0\cr -1/r_k}$ where $r_k$ is the index of  $\Z$ in
$\Z +\Z \mu_{k,k-1}$.
More specifically, one has
$\pmatrix{\delta&-\beta\cr-\gamma&\alpha}\pmatrix{0\cr -1/r_k}=
\pmatrix{1\cr \mu_{k,k-1}}$ so $\beta=r_k$ and
$\alpha=-r_k\mu_{k,k-1}$.
 By \cite{C} we have $|\gamma|\leq |\mu_{k,k-1}r_k|$
and
$|\delta|\leq r_k$.
(Actually this is wrong. Indeed \cite{C} only claims it when $\mu_{k,k-1}$
is nonzero. We leave the modifications for the case $\mu_{k,k-1}=0$
as an exercise.)

Now put $c_{k-1}=\alpha b_{k-1}+ \beta b_k$ and
$c_k= \gamma b_{k-1}+ \delta b_k$. The algorithm \trickle\ tells us to
replace
$b_k$ with $c_k$ and $b_{k-1}$ with $c_{k-1}$.
We want to estimate the resulting new $\mu_{i,j}$, which we call $\nu_{i,j}$.
For $i$ different from $k$, $k-1$ nothing changes.
Further $|\nu_{k-1,j}|= |\alpha\mu_{k-1,j}+\beta\mu_{k,j}|
\leq r_k|\mu_{k,k-1}\mu_{k-1,j}|+r_k|\mu_{k,j}|$ and
$|\nu_{k,j}|= |\gamma\mu_{k-1,j}+\delta\mu_{k,j}|\leq
r_k|\mu_{k,k-1}\mu_{k-1,j}|+r_k|\mu_{k,j}|$, which is the same bound.

\begin{Lemma}
As $k$ descends we have
\begin{enumerate}
\item $|\mu_{i,j}|\leq 1$ for $k> i>j\geq0$,
\item $|\mu_{k,j}|\leq \sqrt B \prod_{i=k+1}^{\kmax+1}(2r_i)$
for $k>j\geq0$,
\item $|\mu_{i,j}|\leq 2^m(\sqrt B)^{\rank+1}$ for $k\leq i>j\geq0$.
\end{enumerate}
\end{Lemma}
\subsubsection*{Proof}
Initially we have $k=\kmax+1$ and $|\mu_{k,j}|^2\leq B$.
Now assume the estimates are true for the present $k$.
We get $|\nu_{k-1,j}|\leq
r_k|\mu_{k,k-1}\mu_{k-1,j}|+r_k|\mu_{k,j}|
\leq 2r_k\max_j|\mu_{k,j}|$ which takes care of $|\nu_{k-1,j}|$.
Now $\prod_{k=\isodim+2}^{\kmax+1}r_k$ is the ratio by which the covolume
drops when adding $e_{\kmax+1}G$ to the lattice spanned by $e_1G$, \dots,
$e_\kmax G$. So it is at most $(\sqrt B)^\rank$.
Thus $|\nu_{k,j}|\leq 2^m(\sqrt B)^{\rank+1}$.
\qed

\subsection{Bailing out}
When $k$ has reached $\isodim+1$
or a new pivot has been created, it is time to forget the old
 $(\ ,\ )_\mix$. But first use the
estimates of the $\mu_{i,j}$ to estimate
$$(b_i,
b_i)_\mix\leq m4^mB^{\rank+1}(\sqrt m(B+1))^{m}$$
and $$(\mu_{i,0}e_{\kmax+1},\mu_{i,0}e_{\kmax+1})_\mix\leq
(B+1)4^mB^{\rank+1},$$
 next
$$(b_i-\mu_{i,0}e_{\kmax+1},
b_i-\mu_{i,0}e_{\kmax+1})\leq m^2(\sqrt m (B+1))^{2m}4^{m+1}B^{\rank+1}
$$ by means of Lemma \ref{mix2euc},
and finally
$$(b_i,b_i)
\leq (4mB)^{4m}$$ say.

Now update $G$, $\isodim$, $\kmax$, $(\ ,\ )_\mix$.
We have to compute the new $\mu_{j,i}$. They can be estimated,
as we have an estimate for $(b_j,b_j)$
 and for
$(b_i^*,b_i^*)_\mix^{-1}$.
We get the estimate $|\mu_{j,i}|^2\leq (4mB)^{5m}$, which was needed in
\ref{preserved}.

\end{document}